\documentclass{amsart}
\usepackage[dvips]{graphicx}
\usepackage{amscd}
\usepackage{amsmath}
\usepackage{amsxtra}
\usepackage{amsfonts}
\usepackage{amssymb}
\newtheorem{theorem}{Theorem}[section]
\newtheorem{corollary}[theorem]{Corollary}
\newtheorem{lemma}[theorem]{Lemma}
\newtheorem{proposition}[theorem]{Proposition}
\theoremstyle{definition}
\newtheorem{definition}[theorem]{Definition}
\newtheorem{remark}[theorem]{Remark}

\theoremstyle{remark}

\renewcommand{\theclaim}{\textup{\theclaim}}

\numberwithin{equation}{section}

\def\openone

{\mathchoice

{\hbox{\upshape \small1\kern-3.3pt\normalsize1}}

{\hbox{\upshape \small1\kern-3.3pt\normalsize1}}

{\hbox{\upshape \tiny1\kern-2.3pt\SMALL1}}

{\hbox{\upshape \Tiny1\kern-2pt\tiny1}}}

\makeatletter

\newbox\ipbox

\newcommand{\ip}[2]{\left\langle #1\,|\,#2\right\rangle}

\newcommand{\diracb}[1]{\left\langle #1\mathrel{\mathchoice

{\setbox\ipbox=\hbox{$\displaystyle \left\langle\mathstrut #1\right.$}

\vrule height\ht\ipbox width0.25pt depth\dp\ipbox}

{\setbox\ipbox=\hbox{$\textstyle \left\langle\mathstrut #1\right.$}

\vrule height\ht\ipbox width0.25pt depth\dp\ipbox}

{\setbox\ipbox=\hbox{$\scriptstyle \left\langle\mathstrut #1\right.$}

\vrule height\ht\ipbox width0.25pt depth\dp\ipbox}

{\setbox\ipbox=\hbox{$\scriptscriptstyle \left\langle\mathstrut #1\right.$}

\vrule height\ht\ipbox width0.25pt depth\dp\ipbox}

}\right. }

\newcommand{\dirack}[1]{\left. \mathrel{\mathchoice

{\setbox\ipbox=\hbox{$\displaystyle \left.\mathstrut #1\right\rangle$}

\vrule height\ht\ipbox width0.25pt depth\dp\ipbox}

{\setbox\ipbox=\hbox{$\textstyle \left.\mathstrut #1\right\rangle$}

\vrule height\ht\ipbox width0.25pt depth\dp\ipbox}

{\setbox\ipbox=\hbox{$\scriptstyle \left.\mathstrut #1\right\rangle$}

\vrule height\ht\ipbox width0.25pt depth\dp\ipbox}

{\setbox\ipbox=\hbox{$\scriptscriptstyle \left.\mathstrut #1\right\rangle$}

\vrule height\ht\ipbox width0.25pt depth\dp\ipbox}

} #1\right\rangle}

\newcommand{\linft}{L^{\infty}\left(  \mathbb{T}\right)}

\newcommand{\loner}{L^{1}\left(\mathbb{R}\right)}

\newcommand{\ltwor}{L^{2}\left(\mathbb{R}\right)}

\newcommand{\ltworn}{L^{2}\left(\mathbb{R}^n\right)}

\newcommand{\ltwozn}{l^{2}\left(\mathbb{Z}^n\right)}
\newcommand{\ltwoz}{l^2\left(\mathbb{Z}\right)}
\newcommand{\Trace}{\operatorname*{Trace}}

\newcommand{\Tper}{\mathcal{T}_{per}}

\newcommand{\Rn}{\mathbb{R}^n}

\newcommand{\cj}[1]{\overline{#1}}

\newcommand{\fo}[1]{\widehat{#1}}

\input cyracc.def

\begin{document}
\title[Local trace for super-wavelets]{The local trace function for super-wavelets}
\author{Dorin Ervin Dutkay}
\address{Department of Mathematics\\
The University of Iowa\\
14 MacLean Hall\\
Iowa City, IA 52242-1419\\
U.S.A.} \email{ddutkay@math.uiowa.edu}
\thanks{}
 \keywords{}

\begin{abstract}
We define an affine structure on $\ltwor\oplus...\oplus\ltwor$
and, following some ideas developed in \cite{Dut1}, we construct a
local trace function for this situation. This trace function is a
complete invariant for a shift invariant subspace and it has a
variety of properties which make it easily computable. The local
trace is then used to give a characterization of super-wavelets
and to analyze their multiplicity function, dimension function and
spectral function. The "$n\times$" oversampling result of Chui and
Shi \cite{CS} is refined to produce super-wavelets.
\end{abstract}

\maketitle \tableofcontents

\section{\label{introduction}Introduction}
The wavelet theory involves the study of an affine structure
existent on a Hilbert space $H$, given by two unitary operators,
$T$ the translation, and $U$ the dilation that satisfy the
commutation relation:
$$UTU^{-1}=T^N,$$
where $N\geq2$ is an integer called the scale. For the classical
wavelet theory, the Hilbert space is $\ltwor$, the translation is
$$T_0f(x)=f(x-1),\quad(x\in\mathbb{R},f\in\ltwor),$$
and the dilation operator is
$$U_0f(x)=\frac{1}{\sqrt{N}}f\left(\frac{x}{N}\right),\quad(x\in\mathbb{R},f\in\ltwor),$$
where $N\geq2$ is an integer.
\par
Consider an affine structure $U$, $T$ on a Hilbert space $H$. A
wavelet is a set $\Psi=\{\psi_1,...,\psi_L\}$ such that the affine
system
$$\{U^mT^k\psi\,|\,m,k\in\mathbb{Z},\psi\in\Psi\}$$ is an
orthonormal basis for $H$. Often the requirement is weaker, and
one looks for the affine system to be only a frame. We recall that
a family $(e_i)_{i\in I}$ of vectors in a Hilbert space $H$ is
called a frame with constants $A>0$ and $B>0$ if
$$A\|f\|^2\leq\sum_{i\in I}|\ip{f}{e_i}|^2\leq B\|f\|^2,\quad(f\in H).$$
When only the second inequality is satisfied, the family is called
Bessel. When $A=B=1$, the frame is called a normalized tight
frame.
\par
Each orthonormal wavelet has a generalized multiresolution
attached to it.
\par
A generalized multiresolution analysis (GMRA) is a collection
$(V_j)_{j\in\mathbb{Z}}$ of subspaces of $H$ that satisfy:
\begin{enumerate}
\item
$V_j\subset V_{j+1}$ for all $j$.
\item
$UV_j=V_{j-1}$ for all $j$.
\item
$\cup V_j$ is dense in $H$ and $\cap V_j=\{0\}$.
\item
$V_0$ is invariant under $T_k$ for all $k\in\mathbb{Z}$.
\end{enumerate}
If, in addition, there is a vector $\varphi\in V_0$ such that
$\{T^k\varphi\,|\,k\in\mathbb{Z}\}$ is an orthonormal basis for
$V_0$, then $(V_j)_j$ is called a multiresolution analysis (MRA)
and $\varphi$ is called a scaling vector for this MRA. For more
information on GMRA's we refer to \cite{BMM}.
\par
The GMRA associated to a wavelet $\Psi$ is defined by
$$V_j:=\{U^mT^k\psi\,|\, m>-j,
k\in\mathbb{Z}\},\quad(j\in\mathbb{Z}).$$ If the GMRA associated
to $\Psi$ is actually a MRA, then $\Psi$ is called a MRA wavelet.
\par
 As it is shown in \cite{HL} and \cite{Dut2}, it is possible to construct
an affine structure and wavelets on the larger Hilbert space
$\ltwor\oplus...\oplus\ltwor$, i.e. "super-wavelets". One way to
do this is by considering the operators:
$$T:=T_0\oplus...\oplus T_0,\quad U:=U_0\oplus...\oplus U_0,$$
and as shown in \cite{HL}, wavelets can be constructed for this
affine structure, but none of them are MRA wavelets.
\par
The affine structure described in \cite{Dut2} takes into
consideration the cycles of the map $z\mapsto z^N$. A set
$\{z_1,...,z_p\}$ of points on the unit circle
$\mathbb{T}:=\{z\in\mathbb{C}\,|\, |z|=1\}$ is called a cycle if
the points are distinct and
$z_1^N=z_2,...,z_{p-1}^N=z_p,z_p^N=z_1$. Then the translation is
defined by
$$T(f_1,...,f_p)=(z_1T_0f_1,...,z_pT_0f_p),\quad(f_1,...,f_p\in\ltwor),$$
and the dilation is defined by
$$U(f_1,...,f_p)=(U_0f_2,U_0f_3,...,U_0f_p,U_0f_1),\quad(f_1,...,f_p\in\ltwor).$$
Note that the dilation operator permutes cyclically the components
$f_1,...,f_p$. Then one can do finite direct sums of these
representations for different cycles to obtain new affine
structures.
\par
The advantage of the second type of affine structure is that it
possesses MRA-wavelets (see \cite{Dut2}).
\par
The objective of this paper is to analyze these affine structures
(see section \ref{perm} for the definition) and their
"super-wavelets", their characterization, their dimension
function, multiplicity function and spectral function, the
relation between super-wavelets and super-scaling functions. We
follow the ideas introduced in \cite{Dut1}, and the main tool will
be the local trace function.
\par
The subspaces of $\ltwor\oplus...\oplus\ltwor$ which are invariant
under all integer translations $T^k$ are called shift invariant.
The local trace function is a complex valued map defined on
$\mathbb{R}$ which is associated to a shift invariant space and an
operator on $\ltwoz\oplus...\oplus\ltwoz$ (definition
\ref{def3_1}). Varying the operator, we obtain a lot of
information about the shift invariant subspace. Taking the
operator to be the identity yields the dimension function
(definition \ref{def2_16}), the spectral function is obtained by
taking the operator to be a certain canonical one-dimensional
projection (definition \ref{def2_12}).
\par
The central results are theorem \ref{th3_3} and \ref{th3_3_1}.
Their power is twofold: they give a formula for computing the
local trace in terms of normalized tight frame generators and also
they show that this formula is independent of the choice of the
normalized tight frame generator. Particular instances of these
theorems give fundamental results about wavelets and shift
invariant subspaces:
\begin{enumerate}
\item
Take the shift invariant subspace to be
$\ltwor\oplus...\oplus\ltwor$ and the consequence is the theorem
that characterizes normalized tight frames generated by
translations: theorem \ref{th2_8}. This in turn gives a
characterization result for super-wavelets, theorem \ref{th2_9}.
\item
Take the shift invariant subspace to be $V_0$, the core space of
the GMRA of a wavelet and take the operator to be $I$, then the
result is the equality between the dimension function and the
multiplicity function (proposition \ref{th2_19}).
\item
For $V_0$, take the operator to be a canonical one-dimensional
projection and the result is the Gripenberg-Weiss formula of
corollary \ref{cor2_15}.
\end{enumerate}
Also, using the spectral function we are able to give a lower
bound estimate on the dimension function (corollary \ref{cor2_22})
which shows not only that some of our affine structures do not
have MRA wavelets, but also we can give a lower bound on the
number of scaling functions any such wavelet needs.
\par
Another interesting application, which follows from the
characterization theorem \ref{th2_9}, is a refinement of an
oversampling result of Chui and Shi \cite{CS}, which states that
when $\psi$ is a NTF wavelet for $N=2$, then
$\eta:=1/p\psi(\cdot/p)$ is also a NTF wavelet, $p$ being an odd
number. We prove in theorem \ref{th2_9_3} that much more is true:
$\eta$ is part of a super-wavelet $\vec\eta:=(\eta,...,\eta)$ in
$\ltwor^p$ and there is also a converse: if $\vec{\eta}$ is a NTF
super-wavelet then $\psi$ is a wavelet. Moreover, going from the
wavelet $\psi$ to the super-wavelet $\vec\eta$ preserves the
orthogonality.
\par
Some notations:
\par
$\mathbb{T}$ is the unit circle in $\mathbb{C}$. We will often
identify it with the interval $[-\pi,\pi)$ and the functions on
$\mathbb{T}$ with functions on $[-\pi,\pi)$ and with
$2\pi$-periodic functions on $\mathbb{R}$. The identification is
done via $z=e^{-i\theta}$.
\par
The Fourier transform is given by
$$\fo{f}(\xi)=\frac{1}{\sqrt{2\pi}}\int_{\mathbb{R}}f(x)e^{-i\xi
x}\,dx,\quad(\xi\in\mathbb{R}).$$ If $T$ in an operator on
$\ltwor$ we denote by $\fo{T}$ its conjugate by the Fourier
transform $\fo{T}\fo{f}=\fo{Tf}$, for $f\in\ltwor$. On
$\ltwor\oplus...\oplus\ltwor$ the Fourier transform is done
componentwise.
\par
For $f\in\linft$, denote by $\pi_0(f)$ the operator on $\ltwor$
defined by
$$\fo{\pi_0(f)}\varphi=f\varphi,\quad(f\in\linft);$$
note that $T_0=\pi_0(z)$ (here $z$ indicates the identity function
on $\mathbb{T}$).

\section{\label{perm}Permutative wavelet representations}
We define here the affine structure on
$\ltwor\oplus...\oplus\ltwor$ that we will work with. This
structure contains the ones used in \cite{HL} and \cite{Dut2}
which we mentioned in the introduction.
\par
Let $\sigma$ be a permutation of the set $\{1,...,n\}$ and
$z_k=e^{-i\theta_k}\in\mathbb{T}$, $\theta_k\in[-\pi,\pi)$, with
the property that $z_k^N=z_{\sigma(k)}$ for all $k\in\{1,...,n\}$.
We denote by $Z:=(z_1,...,z_n)$.
\par
We define a wavelet representation on the Hilbert space
$$\ltwor^n=\underbrace{\ltwor\oplus...\oplus\ltwor}_{n\mbox{
times}}.$$
 For
$f\in\linft$ define the operator $\pi_{\sigma,Z}(f)$ on $\ltwor^n$
by:
$$\pi_{\sigma,Z}(f)(\varphi_1,...,\varphi_n)=(\pi_0(f(zz_1))\varphi_1,...,\pi_0(f(zz_n))\varphi_n),$$
for $\varphi_i\in\ltwor$, ($i\in\{1,...,n\}$), which means that
the Fourier transform of this operator is
$$\widehat{\pi}_{\sigma,Z}(f)(\varphi_1,...,\varphi_n)=(f(\cdot\,+\theta_1)\varphi_1,...,f(\cdot\,+\theta_n)\varphi_n).$$
Define the dilation $U_{\sigma,Z}$ on $\ltwor^n$ by
$$U_{\sigma,Z}(\varphi_1,...,\varphi_n)=(U_0\varphi_{\sigma(1)},...,U_0\varphi_{\sigma(n)}).$$
Then $(\ltwor^n,\pi_{\sigma,Z},U_{\sigma,Z})$ is a normal wavelet
representation, i.e. it satisfies the following conditions:
\begin{enumerate}
\item
$U_{\sigma,Z}$ is unitary;
\item
$\pi_{\sigma,Z}$ is a unital representation of the $C^*$-algebra
$\linft$;
\item
$U_{\sigma,Z}\pi_{\sigma,Z}(f)U_{\sigma,Z}^{-1}=\pi_{\sigma,Z}(f(z^N))$
for all $f\in\linft$;
\item
For every uniformly bounded sequence $(f_n)_{n\in\mathbb{N}}$ in
$\linft$ which converges a.e to a function $f\in\linft$, one has
that $\pi_{\sigma,Z}(f_n)$ converges to $\pi_{\sigma,Z}(f)$ in the
strong operator topology.
\end{enumerate}
Checking these properties requires just some simple computations.
The translation (or shift) of this wavelet representation is
$$T_{\sigma,Z}:=\pi_{\sigma,Z}(z),$$
where $z$ indicates the identity function on $\mathbb{T}$,
$z\mapsto z$. Also note that $T_{\sigma,Z}^n=\pi_{\sigma,Z}(z^n)$.
\begin{remark}\label{rem1_1}
Each permutation $\sigma$ can be decomposed into a finite product
of disjoint cycles. Since $z_i^N=z_{\sigma(i)}$ for all $i$, this
implies that all $z_i$ have a finite orbit under the map $z\mapsto
z^N$, i.e. they are cycles, as defined in the introduction.
\par
It may be that the length of the cycle given by $z_i$ is shorter
then the length of the corresponding cycle of $\sigma$ (for
example, when all $z_i=1$ and $\sigma$ is a product of disjoint
transpositions, then the length of the cycle of $z_i$ is $1$, but
the cycles of $\sigma$ have length $2$). However when these
lengths coincide we obtain the representations defined in
\cite{Dut2}.
\par
When $z_i=1$ for all $i$ and $\sigma$ is the identity we have the
amplification of the standard representation on $\ltwor$,
amplification which was studied in \cite{HL}.
\end{remark}
\par
For the rest of the paper we will consider a fixed permutation
$\sigma$ and the points $z_1,z_2,...,z_n$, and we will omit the
subscripts.

\section{\label{shif}{Shift invariant subspaces}}
In this section we present some structural theorem for subspaces
of $\ltwor^n$ which are invariant under the integer translations
$T^k$.
\begin{definition}
A subspace $V$ of $\ltwor^n$ is called shift invariant (SI) if
$T^k V\subset V$ for all $k\in\mathbb{Z}$.
\par
If $\mathcal{A}$ is a subset of $\ltwor^n$, we denote by
$S(\mathcal{A})$ the shift invariant subspace generated by
$\mathcal{A}$.
\end{definition}
\begin{definition}
Let $V$ be a shift invariant subspace of $\ltworn$. A subset
$\Phi$ of $V$ is called a normalized tight frame generator (or NTF
generator) for $V$ if
$$\{T_k\varphi\,|\,k\in\mathbb{Z},\varphi\in\Phi\}$$
is a NTF for $V$.
\par
We use also the notation $S(\varphi):=S(\{\varphi\})$. $\varphi$
is called a quasi-orthogonal generator for $S(\varphi)$ if
$$\{T_k\varphi\,|\,k\in\mathbb{Z}\}$$
is a NTF for $S(\varphi)$ and for all $\xi\in\mathbb{R}$,
$$\operatorname*{Per}|\widehat{\varphi}|^2(\xi):=\sum_{k\in\mathbb{Z}}|\widehat{\varphi}|^2(\xi+2k\pi)\in\{0,1\}.$$
(actually, the second condition is a consequence of the first but
we include it anyway).
\end{definition}

\par
We will need the "fiberization" techniques based on the range
function of Helson \cite{H}, which were used in the classical case
in the work of de Boor, DeVore and Ron \cite{BDR} and further
developed by Ron and Shen \cite{RS1}-\cite{RS4}, Bownik \cite{Bo1}
and others.
\par
For $\varphi=(\varphi_1,...,\varphi_n)\in\ltwor^n$, the map
$\Tper\varphi$ assigns to each point $\xi\in\mathbb{R}$ a vector
in $\ltwoz^n$, called the fiber of $\varphi$ at $\xi$:
\par
$$\Tper(\varphi_1,...,\varphi_n)(\xi)=\left((\widehat{\varphi}_1(\xi-\theta_1+2k\pi))_{k\in\mathbb{Z}},...,(\widehat{\varphi}_n(\xi-\theta_n+2k\pi))_{k\in\mathbb{Z}}\right),\,(\xi\in\mathbb{R}),$$
for all $(\varphi_1,...,\varphi_n)\in\ltwor^n$.
\par
This map transforms the translations into multiplications by
scalars in each fiber:
\begin{proposition}
For all $f\in\linft$ and $(\varphi_1,...,\varphi_n)\in\mathbb{R}$,
$$\Tper(\pi(f)(\varphi_1,...,\varphi_n))(\xi)=f(\xi)\Tper(\varphi_1,...,\varphi_n)(\xi),\quad(\xi\in\mathbb{R}).$$
In particular, for all $k\in\mathbb{Z}$,
$$\Tper(T^k(\varphi_1,...,\varphi_n)(\xi)=e^{-ik\xi}\Tper(\varphi_1,...,\varphi_n)(\xi),\quad(\xi\in\mathbb{R}).$$
\end{proposition}
\par
Note also the periodicity property of $\Tper$ (which justifies the
subscript "per"):
$$\Tper\varphi(\xi+2s\pi)=\lambda(s)^*(\Tper f(\xi)),\quad(\xi\in\mathbb{R},s\in\mathbb{Z}),$$
where, for $s\in\mathbb{Z}$, $\lambda(s)$ is the shift on
$\ltwoz^n$:
$$(\lambda(s)\alpha)(k,i)=\alpha(k-s,i),\quad(k\in\mathbb{Z},i\in\{1,...,N\}).$$
\par
For a shift invariant subspace the fibers $\Tper\varphi(\xi)$ at a
fixed point $\xi$ form a subspace of $\ltwoz^n$, which will be
denoted by $J_{per}(\xi)$. Thus $J_{per}$ will be a bundle map
which assigns to each point $\xi\in\mathbb{R}$ a subspace of
$\ltwoz^n$. This is the idea behind the range function.
\begin{definition}
A range function is a measurable mapping
$$J:[-\pi,\pi)\rightarrow\{\mbox{ closed subspaces of }\ltwoz^n\}.$$
Measurable means weakly operator measurable, i.e.,
$\xi\mapsto\ip{P_{J(\xi)}a}{b}$ is measurable for any choice of
vectors $a,b\in\ltwoz^n$.
\par
A periodic range function is a measurable function
$$J_{per}:\mathbb{R}\rightarrow\{\mbox{ closed subspaces of }\ltwoz^n\},$$
with the periodicity property:
$$J_{per}(\xi+2k\pi)=\lambda(k)^*\left(J_{per}(\xi)\right),\quad(k\in\mathbb{Z},\xi\in\mathbb{R}).$$
\end{definition}
 Sometimes we will use the same letter to denote the subspace $J_{per}(\xi)$ and the projection onto
$J_{per}(\xi)$. In terms of projections, the periodicity can be
written as:
$$J_{per}(\xi+2k\pi)=\lambda(k)^*J_{per}(\xi)\lambda(k),\quad(k\in\mathbb{Z},\xi\in\mathbb{R}).$$
\par
There is a one-to-one correspondence between shift invariant
subspaces and range functions. This correspondence is made precise
in the next theorem due to Helson \cite{H}. The proof given in
\cite{Bo1}, proposition 1.5, carries over here, without any
significant modification.
\begin{theorem}
A closed subspace $V$ of $\ltworn$ is shift invariant if and only
if
$$V=\{f\in\ltwor^n\,|\, \Tper f(\xi)\in J_{per}(\xi)\mbox{ for a.e. }\xi\in\mathbb{R}^n\},$$
for some measurable periodic range function $J_{per}$. The
correspondence between $V$ and $J_{per}$ is bijective under the
convention that range functions are identified if they are equal
a.e. Furthermore, if $V=S(\mathcal{A})$ for some countable
$\mathcal{A}\subset\ltworn$, then
$$J_{per}(\xi)=\overline{\operatorname*{span}}\{\Tper\varphi(\xi)\,|\,\varphi\in\mathcal{A}\},\quad \mbox{for
a.e. }\xi\in\mathbb{R}.$$
\end{theorem}
\par
Another fundamental fact about the range function is that it
transforms the frame property into a local property: more
precisely, a family of vectors generate a normalized tight frame
by translations if and only if the fibers form normalized tight
frames at each point in $\mathbb{R}$. This is given in the next
theorem, for a proof look in \cite{Bo1}, theorem 2.3.
\begin{theorem}\label{th2_5}
Let $V$ be a SI subspace of $\ltwor^n$, $J_{per}$ its periodic
range function and $\Phi$ a countable subset of $V$.
$\{T^k\varphi\,|\,k\in\mathbb{Z},\varphi\in\Phi\}$ is a frame with
constants $A$ and $B$ for $V$ (Bessel family with constant $B$) if
and only if $\{\Tper\varphi(\xi)\,|\,\varphi\in\Phi\}$ is a frame
with constants $A$ and $B$ for $J_{per}(\xi)$ (Bessel sequence
with constant $B$) for almost every $\xi\in\mathbb{R}$.
\end{theorem}
\par
Another approach, with a more representation-theoretic flavor, is
proposed in \cite{BMM},\cite{BM},\cite{B}. The group $\mathbb{Z}$
has a unitary representation on a shift invariant space and the
abstract harmonic analysis arguments given in the papers by
Baggett et al work here with minor changes and yield the following
results.
\begin{theorem}\label{th2_6}
If $V_0$ is a shift invariant subspace of $\ltwor^n$ then there
exists a multiplicity function
$m:\mathbb{T}\rightarrow\{0,1,...,\infty\}$ such that:
\begin{enumerate}
\item
If
$$S_j:=\{\xi\in\mathbb{T}\,|\,m(\xi)\geq j\},$$
then there is a unitary
$$\mathcal{J}:V_0\rightarrow\bigoplus_{j=1}^\infty L^2(S_j)$$
which intertwines the representations
\begin{equation}\label{eq2_6_1}
\mathcal{J}T^k=\rho_k\mathcal{J},\quad(k\in\mathbb{Z}),
\end{equation}
where $\rho_k$ is multiplication by the function $\xi\mapsto
e^{-ik\xi}$ on each component:
$$\rho_k(f_1,f_2,...)=(e^{-ik\cdot}f_1,e^{-ik\cdot}f_2,...),\quad(k\in\mathbb{Z}).$$
\item
If
$$\phi_i:=\mathcal{J}^{-1}(\chi_{S_i}),\quad(i\in\{1,2...\}),$$
(where $\chi_{S_i}$ stands for the element of $\oplus L^2(S_j)$
whose only nonzero component is $\chi_{S_i}$ in $L^2(S_i)$), then
\begin{equation}\label{eq2_6_2}
\{\phi_i\,|\,i\in\{1,2,...\}\}\mbox{ is a NTF generator for }V_0,
\end{equation}
\begin{equation}\label{eq2_6_3}
T^k\phi_i\perp T_l\phi_j\quad(i\neq j,k,l\in\mathbb{Z}),
\end{equation}
and
\begin{equation}\label{eq2_6_4}
\ip{\Tper\phi_i(\xi)}{\Tper\phi_j(\xi)}=\left\{\begin{array}{ccc}
\chi_{S_i},&\mbox{if}&i=j\\
0,&\mbox{if}&i\neq j.
\end{array}\right.\quad(\xi\in[-\pi,\pi).
\end{equation}
\item
$V_0$ can be decomposed as the orthogonal sum of the shift
invariant subspaces generated by the quasi-orthogonal generators
$\phi_i$:
\begin{equation}\label{eq2_6_5}
V_0=\bigoplus_{i=1}^\infty S(\phi_i).
\end{equation}
\end{enumerate}
\end{theorem}
\begin{proof}
We follow the arguments given in \cite{BMM} which use the spectral
theory of Stone and Mackey for the commutative group $\mathbb{Z}$
that has as dual $\widehat{\mathbb{Z}}=\mathbb{T}$.
\par
For the representation of $\mathbb{Z}$ by translations $T^k$ on
$V_0$, there is a unique projection-valued measure $p$ on
$\mathbb{T}$ for which
$$T^k|_{V_0}=\int_{\mathbb{T}}e^{-ik\xi}\,dp(\xi),\quad(k\in\mathbb{Z}).$$
\par
The projection-valued measure $p$ is completely determined by a
measure class $[\mu]$ on $\mathbb{T}$ and a measurable
multiplicity function $m:\mathbb{T}\rightarrow\{0,1,...,\infty\}$.
\par
As in \cite{BMM} proposition 1.2, we will show that the measure
class $[\mu]$ is absolutely continuous with respect to the Haar
measure on $\mathbb{T}$.
\par
Indeed, the Fourier transform establishes an equivalence between
the translation $T$ on $\ltwor^n$ and its Fourier version
$$\widehat{T}^k(f_1,...,f_n)=(e^{-ik(\cdot\,+\theta_i)}f_i)_{i\in\{1,...,n\}},\quad(k\in\mathbb{Z}).$$
Define $W:\ltwor^n\rightarrow\ltwor^n$ by
$$W(f_1,...,f_n)=(f_i(\cdot-\theta_i))_{i\in\{1,...,n\}}.$$
This intertwines the given representation with the one defined by
$$\tilde{T}^k(f_1,...,f_n)=(e^{-ik\cdot }f_i)_{i\in\{1,...,n\}}.$$
But $\mathbb{R}$ can be decomposed as
$$\mathbb{R}=\cup_{k\in\mathbb{Z}}[(2k-1)\pi,(2k+1)\pi),$$
thus, the representation of $\mathbb{Z}$ on the entire space
$\ltwor^n$ is equivalent to $n$ times an infinite multiple of the
regular representation of $\mathbb{Z}$ (which is multiplication by
$e^{-ik\xi}$ on $L^2(\mathbb{T})$.
\par
Therefore, the projection-valued measure associated to the
representation of $\mathbb{Z}$ on the entire space $\ltwor^n$ is
equivalent to the Haar measure on $\mathbb{T}$, so the
projection-valued measure associated to any subrepresentation of
this representation must have a measure which is absolutely
continuous w.r.t. the Haar measure.
\par
For (\ref{eq2_6_1}), (\ref{eq2_6_2}) and (\ref{eq2_6_3}) the
argument is the one used on \cite{BM} lemma 1.1.
\par
(\ref{eq2_6_4}) can be proved as in lemma 1.2 from \cite{BM}. It
follows from the fact that $\mathcal{J}$ is unitary and
intertwining so, for $k\in\mathbb{Z},i,j\in\{1,2,...\}$,
$$\ip{T^k\phi_i}{\phi_j}=\int_{\mathbb{T}}e^{-ik\xi}\chi_{S_i}\delta_{ij}\,d\xi,$$
which implies
$$\int_{\mathbb{T}}e^{-ik\xi}\ip{\Tper\phi_i(\xi)}{\Tper\phi_j(\xi)}d\xi=\int_{\mathbb{T}}e^{-ik\xi}\chi_{S_i}\delta_{ij}\,d\xi.$$
\par
(\ref{eq2_6_5}) is just a rephrasing of the previous statements.
\end{proof}

\section{\label{loca}The local trace function}
We are able now to define the local trace function. The proofs
from \cite{Dut1} work here in most of the instances and, when a
result is presented without a proof, the reader can look in
\cite{Dut1}.
\par
Each shift invariant subspace $V$, can be seen as a bundle of
vector subspaces of $\ltwoz^n$ by means of a range function
$J_{per}$. The local trace function associated to $V$ and an
operator $T$ on $\ltwoz^n$ is the function which associates to
each point $\xi$ the value of the canonical trace of the
restriction of $T$ to the subspace $J_{per}(\xi)$. We recall that
the trace of a positive or trace-class operator $T$ on some
Hilbert space $H$, can be computed by
$$\Trace(T)=\sum_{i\in I}\ip{Tf_i}{f_i},$$
where $(f_i)_{i\in I}$ is any normalized tight frame of the
Hilbert space $H$ (see \cite{Dut1}).
\begin{definition}\label{def3_1}
Let $V$ be a SI subspace of $\ltworn$, $T$ a positive (or
trace-class) operator on $\ltwoz^n$ and let $J_{per}$ be the range
function associated to $V$. We define the local trace function
associated to $V$ and $T$ as the map from $\mathbb{R}$ to
$[0,\infty]$ (or $\mathbb{C}$) given by the formula
$$\tau_{V,T}(\xi)=\Trace\left(TJ_{per}(\xi)\right),\quad(\xi\in\mathbb{R}).$$
We define the restricted local trace function associated to $V$
and a vector $f$ in $\ltwoz^n$ by
$$\tau_{V,f}(\xi)=\Trace\left(P_fJ_{per}(\xi)\right)(=\tau_{V,P_f}(\xi)),\quad(\xi\in\mathbb{R}^n),$$
where $P_f$ is the operator on $\ltwoz^n$ defined by
$P_f(v)=\ip{v}{f}f$.
\end{definition}
\begin{remark}
Even though, $TJ_{per}(\xi)$ is not necessarily positive when $T$
is, when computing the trace, we can consider instead of
$TJ_{per}(\xi)$ the operator\\ $J_{per}(\xi)TJ_{per}(\xi)$, or we
can use normalized tight frames just for the subspace
$J_{per}(\xi)$ (not for the entire $\ltwoz^n$), in any case the
expression $\Trace(TJ_{per}(\xi))$ makes sense and is a positive
number.
\par
When $T$ is trace-class, since the trace class operators form an
ideal, $TJ_{per}(\xi)$ is also trace class so
$\Trace(TJ_{per}(\xi))$
\end{remark}
\par
Here are some ways to compute the local trace function.
\begin{proposition}\label{prop3_1}
For all $f\in\ltwozn$,
$$\tau_{V,f}(\xi)=\|J_{per}(\xi)(f)\|^2,\quad(\mbox{ for a.e. }\xi\in\mathbb{R}).$$
\end{proposition}
\par
The next result is central for the paper. As we have mentioned in
the introduction, this theorem can be used for two purposes: it
gives a method to compute the local trace function and secondly,
and maybe even more important, it shows that the formula on the
right does not depend on the choice of the normalized tight frame.
So one can look at the local trace function as an invariant for
the SI subspace, and the theorem is an index theorem.
\begin{theorem}\label{th3_3}
Let $V$ be a SI subspace of $\ltwor^n$ and $\Phi\subset V$ a NTF
generator for $V$. Then for every positive (or trace-class)
operator $T$ on $\ltwoz^n$ and any $f\in\ltwoz^n$,
\begin{equation}\label{eq3_3_1}
\tau_{V,T}(\xi)=\sum_{\varphi\in\Phi}\ip{T\Tper\varphi(\xi)}{\Tper\varphi(\xi)},\quad(\mbox{
for a.e. }\xi\in\mathbb{R});
\end{equation}
\begin{equation}\label{eq3_3_2}
\tau_{V,f}(\xi)=\sum_{\varphi\in\Phi}|\ip{f}{\Tper\varphi(\xi)}|^2,\quad
\mbox{ for a.e. }\xi\in\mathbb{R}).
\end{equation}
\end{theorem}
The converse is also true: if a family of vectors can be used to
compute the local trace then they form a NTF generator for the SI
space. Moreover it is sufficient that this is satisfied only for
some very particular vectors (see (iii) in the next result):
\begin{theorem}\label{th3_3_1}
Let $V$ be a SI subspace of $\ltwor^n$, $J_{per}$ its periodic
range function and $\Phi$ a countable subset of $\ltworn$. Then
following affirmations are equivalent:
\begin{enumerate}
\item
$\phi\subset V$ and $\phi$ is a NTF generator for $V$;
\item
For every $f\in\ltwoz^n$
\begin{equation}\label{eq3_3_1_1}
\sum_{\varphi\in\Phi}|\ip{f}{\Tper\varphi(\xi)}|^2=\|J_{per}(\xi)(f)\|^2,\quad\mbox{for
a.e. }\xi\in\mathbb{R}
\end{equation}
\item
For every $0\neq k\in\mathbb{Z}$, $l,j\in\{1,...,n\}$ and
$\alpha\in\{0,1,i=\sqrt{-1}\}$,
\begin{equation}\label{eq3_3_1_2}
\sum_{\varphi\in\Phi}|\widehat{\varphi}_l(\xi)+\overline{\alpha}\widehat{\varphi}_j(\xi+2k\pi)|^2=\|J_{per}(\xi)(\delta_{0l}+\alpha\delta_{kj})\|^2,\quad\mbox{for
a.e. }\xi\in\mathbb{R}.
\end{equation}
\end{enumerate}
\end{theorem}
The local trace function is related in a very simple way to
another important invariant, the dual Gramian. This was introduced
by A. Ron and Z. Shen and successfully used for the analysis of
the structure of SI spaces and of frames generated by translations
in \cite{RS1}-\cite{RS4}.
\begin{definition}\label{def3_4}
Let $V$ be a shift invariant subspace of $\ltwor^n$ and $\Phi$, a
NTF generator for $V$. The dual Gramian of $\Phi$ is the function
$\tilde{G}^{\Phi}$ which assigns to each point $\xi\in\mathbb{R}$
the matrix
$$\tilde{G}^{\Phi}_{ki,lj}(\xi)=\sum_{\varphi\in\Phi}\fo{\varphi}_i(\xi-\theta_i+2k\pi)\cj{\fo{\varphi}}_j(\xi-\theta_j+2l\pi),\quad(k,l\in\mathbb{Z},i,j\in\{1,...,n\}).$$
\end{definition}
We can recuperate the Gramian from the local trace function just
by computing it at some rank-one operators.
\begin{proposition}\label{prop3_4_1}
If $V$ is a SI space, $J_{per}$ its range function, and $\Phi$ a
NTF generator for it then, for all
$k,l\in\mathbb{Z},i,j\in\{1,...,n\}$,
$$\tilde{G}^{\Phi}_{ki,lj}(\xi)=\tau_{V,P_{ki,lj}}(\xi)=\ip{J_{per}(\xi)\delta_{lj}}{J_{per}(\xi)\delta_{ki}},\quad(\xi\in\mathbb{R}),$$
where $P_{ki,lj}$ is the rank-one operator defined by
$$P_{ki,lj}v=\ip{v}{\delta_{ki}}\delta_{lj},\quad(v\in\ltwoz^n).$$
\end{proposition}
\begin{proof}
The first part is a consequence of theorem \ref{th3_3}. For the
last equality we use the following argument: if $(f_q)_{q\in Q}$
is an orthonormal basis for $J_{per}(\xi)$ then
\begin{align*}
\Trace(P_{ki,lj}J_{per}(\xi))&=\sum_{q\in
Q}\ip{P_{ki,lj}f_q}{f_q}\\
&=\sum_{q\in Q}\ip{f_q}{\delta_{ki}}\ip{\delta_{lj}}{f_q}\\
&=\sum_{q\in Q}\ip{f_q}{J_{per}(\xi)\delta_{ki}}\ip{J_{per}(\xi)\delta_{lj}}{f_q}\\
&=\ip{J_{per}(\xi)\delta_{lj}}{J_{per}(\xi)\delta_{ki}}.
\end{align*}
\end{proof}
\par
Having these it is clear that the dual Gramian is an invariant and
in fact it is a complete one (see \cite{Dut1},theorem 4.6) for the
proof:
\begin{theorem}\label{th3_3_2}
Let $V$ be a SI subspace of $\ltwor^n$, $\Phi_1$ a NTF generator
for $V$ and $\Phi_2$ a countable family of vectors from
$\ltwor^n$. The following affirmations are equivalent:
\begin{enumerate}
\item
$\Phi_2\subset V$ and $\Phi_2$ is a NTF generator for $V$;
\item
The dual Gramians $\tilde{G}^{\Phi_1}$ and $\tilde{G}^{\Phi_2}$
are equal almost everywhere.
\end{enumerate}
\end{theorem}
\par
If we particularize this result to the case when the shift
invariant subspace is $\ltwor^n$, then since the range function
for it is constant $\ltwoz^n$ at each point, we obtain the
following characterization of families that generate normalized
tight frames by translations:
\begin{theorem}\label{th2_8}
Let $\Phi$ be a countable subset of $\ltwor^n$. The following
affirmations are equivalent:
\begin{enumerate}
\item
$\Phi$ is a NTF generator for $\ltwor^n$;
\item
$$\sum_{\varphi\in\Phi}\fo{\varphi}_i(\xi-\theta_i)\cj{\fo{\varphi}}_j(\xi-\theta_j+2k\pi)=\delta_{ij}\delta_k,\quad(\xi\in\mathbb{R},i\in\{1,...,n\},k\in\mathbb{Z}).$$
\end{enumerate}
\end{theorem}

\par
Next, we present some elementary properties of the local trace
function. The proofs of these results are simple and can be found
in \cite{Dut1} for the classical case.
\begin{proposition}\label{prop3_4}
{\bf [Periodicity]} Let $V$ be a SI subspace, $T$ a positive or
trace-class operator on $\ltwoz^n$, $f\in\ltwoz^n$. Then, for
$k\in\mathbb{Z}$
\begin{equation}\label{eq3_4_1}
\tau_{V,T}(\xi+2k\pi)=\tau_{V,\lambda(k)T\lambda(k)^*}(\xi),\quad
(\mbox{ for a.e. }\xi\in\mathbb{R});
\end{equation}
\begin{equation}\label{eq3_4_2}
\tau_{V,f}(\xi+2k\pi)=\tau_{V,\lambda(k)f}(\xi),\quad (\mbox{ for
a.e. }\in\mathbb{R}),
\end{equation}
where $\lambda(k)$ is the shift on $\ltwoz^n$, i.e. the unitary
operator on $\ltwoz^n$ defined by
$(\lambda(k)\alpha)(l,i)=\alpha(l-k,i)$ for all
$l\in\mathbb{Z},i\in\{1,...,n\}$, $\alpha\in\ltwoz^n$.
\end{proposition}

\begin{proposition}\label{prop3_6}
{\bf [Additivity]} Suppose $(V_i)_{i\in I}$ are mutually
orthogonal SI subspaces ($I$ countable) and let $V=\oplus_{i\in
I}V_i$. Then, for every positive or trace-class operator $T$ on
$\ltwoz^n$ and every $f\in\ltwoz^n$:
\begin{equation}\label{eq3_6_1}
\tau_{V,T}=\sum_{i\in I}\tau_{V_i,T},\quad\mbox{a.e. on
}\mathbb{R};
\end{equation}
\begin{equation}\label{eq3_6_2}
\tau_{V,f}=\sum_{i\in I}\tau_{V_i,f},\quad\mbox{a.e. on
}\mathbb{R}.
\end{equation}
\end{proposition}
\begin{proposition}\label{prop3_7}
{\bf [Monotony and injectivity]} Let $V,W$ be SI subspaces.
\par
(i) $V\subset W$ iff $\tau_{V,T}\leq\tau_{W,T}$ a.e. for all
positive operators $T$ iff $\tau_{V,f}\leq\tau_{W,f}$ a.e. for all
$f\in\ltwoz^n$.
\par
(ii) $V=W$ iff $\tau_{V,f}=\tau_{W,f}$ a.e for all $f\in\ltwoz^n$.
\end{proposition}

\begin{theorem}\label{th4_3}
Let $(V_j)_{j\in\mathbb{N}}$ be an increasing sequence of SI
subspaces of $\ltworn$,
$$V:=\overline{\bigcup_{j\in\mathbb{N}}V_j},$$
$T$ a positive (trace-class) operator on $\ltwoz^n$ and
$f\in\ltwoz^n$. Then, for a.e. $\xi\in\mathbb{R}$,
$\tau_{V_j,T}(\xi)$ increases (converges) to $\tau_{V,T}(\xi)$ and
$\tau_{V_j,f}(\xi)$ increases to $\tau_{V,f}(\xi)$.
\end{theorem}

\par
We will use the local trace function for the analysis of wavelets.
For this reason we have to see how it behaves under dilation.
\begin{theorem}\label{th2_7}
Let $V$ be a shift invariant subspace of $\ltwor^n$ and $T$ a
positive or trace-class operator on $\ltwoz^n$. For
$l\in\{0,...,N-1\}$ define the operators $D_l$ on $\ltwoz^n$ by
$$D_lv(k,i)=\left\{\begin{array}{ccc}
v(s,i),&\mbox{if}&2k\pi=2l\pi+\theta_i-N\theta_{\sigma^{-1}(i)}+2Ns\pi,\mbox{
with }s\in\mathbb{Z},\\
0,& &\mbox{otherwise},
\end{array}\right.
$$
and the unitary operator $\mathcal{S}$ on $\ltwoz^n$, by
$$\mathcal{S}(v_1,...,v_n)=(v_{\sigma^{-1}(1)},...,v_{\sigma^{-1}(n)}),\quad
(v_1,...,v_n)\in\ltwoz^n).$$ Then $U^{-1}V$ is shift invariant and
its local trace function is given by
$$\tau_{U^{-1}V,T}(\xi)=\sum_{l=0}^{N-1}\tau_{V,\mathcal{S}^*D_l^*TD_l\mathcal{S}}\left(\frac{\xi+2l\pi}{N}\right),\quad(\xi\in\mathbb{R}).$$
\end{theorem}

\begin{proof}
By theorem \ref{th2_6}, $V$ can be decomposed into an orthogonal
sum of singly-generated SI spaces. Also the local trace is
additive and $U$ is unitary so we may assume that $V=S(\varphi)$
with $\varphi$ a quasi-orthogonal generator.
\par
Due to the commutation relation $UTU^{-1}=T^N$, the vectors
$$\{\varphi_{-1,r}:=U^{-1}T^r\varphi,|\,r\in\{0,...,N-1\}\}$$ form a NTF
generator for $U^{-1}V$. Their Fourier transforms are
$$\widehat{\varphi}_{-1,r}(\xi)=\left(\frac{1}{\sqrt{N}}e^{-ir(\xi/N+\theta_{\sigma^{-1}(i)})}\widehat{\varphi}_{\sigma^{-1}(i)}(\xi/N)\right)_{i\in\{1,...,n\}},\quad(\xi\in\mathbb{R}).$$
Then a short computation shows that, if we consider the vectors
$v_l(\xi)$, $l\in\{0,...,N-1\}$ in $\ltwoz^n$,
$$v_l(\xi)(k,i)=\left\{\begin{array}{ccc}
\widehat{\varphi}_{\sigma^{-1}(i)}\left(\frac{\xi-\theta_i+2k\pi}{N}\right),&\mbox{if}&\frac{-\theta_i+2k\pi}{N}\in\frac{2l\pi}{N}-\theta_{\sigma^{-1}(i)}+2\pi\mathbb{Z},\\
0,& &\mbox{otherwise},
\end{array}
\right.$$
then
$$\Tper\varphi_{-1,r}(\xi)=\frac{1}{\sqrt{N}}e^{-ir\frac{\xi}{N}}\sum_{l=0}^{N-1}e^{-ir\frac{2\pi l}{N}}v_l(\xi),\quad(r\in\{0,...,N-1\}).$$
The matrix $$\left(\frac{1}{\sqrt{N}}e^{-ir\frac{2\pi
l}{N}}\right)_{l,r\in\{0,...,N-1\}}$$ is unitary, so the vectors
$\{v_l(\xi)\,|\,l\in\{0,...,N-1\}\}$ span the same subspace of
$\ltwoz^n$ as the vectors
$\{\Tper\varphi_{-1,r}(\xi)\,|\,r\in\{0,...,N-1\}\}$, namely
$J_{U^{-1}V}(\xi)$.
\par
Note that
$$v_l(\xi)=D_l\mathcal{S}\Tper\varphi\left(\frac{\xi+2l\pi}{N}\right),\quad(l\in\{0,...,N-1\}),$$
$D_l$ are isometries and $\mathcal{S}$ is unitary. In addition
$\varphi$ is quasi-orthogonal so $\|\Tper\varphi(\xi)\|\in\{0,1\}$
for almost every $\xi\in\mathbb{R}$, hence
$\|v_l(\xi)\|\in\{0,1\}$ for all $l$. It is also clear that the
vectors $v_l(\xi)$ are mutually orthogonal for a fixed $\xi$.
Having all these, it follows that
$\{v_l(\xi)\,|\,l\in\{0,...,N-1\}\}$ is a NTF for
$J_{U^{-1}V}(\xi)$, therefore we can compute the trace with them.
\begin{align*}
\tau_{U^{-1}V,T}(\xi)&=\sum_{l=0}^{N-1}\ip{Tv_l(\xi)}{v_l(\xi)}\\
&=\sum_{l=0}^{N-1}\ip{TD_l\mathcal{S}\Tper\varphi\left(\frac{\xi+2l\pi}{N}\right)}{D_l\mathcal{S}\Tper\varphi\left(\frac{\xi+2l\pi}{N}\right)}\\
&=\sum_{l=0}^{N-1}\ip{\mathcal{S}^*D_l^*TD_l\mathcal{S}\Tper\varphi\left(\frac{\xi+2l\pi}{N}\right)}{\Tper\varphi\left(\frac{\xi+2l\pi}{N}\right)},
\end{align*}
which proves the formula.
\end{proof}

\section{\label{char}A characterization of super-wavelets}
We apply now the local trace function to super-wavelets. Our first
goal is to obtain a characterization of NTF super-wavelets. This
is done in theorem \ref{th2_9}.
\begin{definition}
A finite subset $\Psi=\{\psi^1,...,\psi^L\}$ of $\ltwor^n$ is a
NTF wavelet if the affine system
$$X(\Psi):=\{\psi_{j,k}:=U^{-j}T^k\psi\,|\,j\in\mathbb{Z},k\in\mathbb{Z},\psi\in\Psi\}$$
is a NTF for $\ltwor^n$.
\par
The quasi-affine system $X^q(\Psi)$ is
$$X^q(\Psi):=\{\tilde{\psi}_{j,k}\,|j\in\mathbb{Z},k\in\mathbb{Z},\psi\in\Psi\},$$
with the convention
$$\tilde{\psi}_{j,k}:=\left\{\begin{array}{ccc}
U^{-j}T^k\psi,&\mbox{if}&j\geq0,k\in\mathbb{Z}\\
N^{j/2}T^kU^{-j}\psi&\mbox{if}&j<0,k\in\mathbb{Z}.
\end{array}\right.$$
\end{definition}
There is an equivalence between affine and quasi-affine frames.
This was proved in full generality for the $\ltwor$-case in
\cite{CSS}. Again the result will work here with some slight
modifications in the proof.
\begin{theorem}\label{th2_9_0}
Let $\Psi=\{\psi^1,...,\psi^L\}$ be a finite subset of $\ltwor^n$.
\begin{enumerate}
\item
$X(\Psi)$ is a Bessel family if and only if $X^q(\Psi)$ is a
Bessel family. Furthermore, their exact frame bounds are equal.
\item
$X(\Psi)$ is an affine frame with constants $A$ and $B$ if and
only if $X^q(\Psi)$ is a quasi-affine frame with constants $A$ and
$B$.
\end{enumerate}
\end{theorem}
\begin{proof}
The proof of theorem 2 in \cite{CSS} works here word for word once
we have the following lemma, which is the analogue of lemma 4 in
\cite{CSS}:
\begin{lemma}
Let $\Psi=\{\psi_1,...,\psi_L\}\subset\ltwor^n$ and let
$f\in\ltwor^n$ be a vector such that $f_1,...,f_n$ have all
compact supports. Then
$$\lim_{M\rightarrow\infty}\sum_{j<0}\sum_{\psi\in\Psi}\sum_{k\in\mathbb{Z}}|\ip{U^{-M}f}{\tilde{\psi}_{j,k}}|^2=0,$$
$$\lim_{M\rightarrow\infty}N^{-M}\sum_{j\leq
-M}\sum_{\nu=0}^{N^M-1}\sum_{\psi\in\Psi}\sum_{k\in\mathbb{Z}}|\ip{T^\nu
f}{\psi_{j,k}}|^2=0.$$
\end{lemma}

Lemma 4 from \cite{CSS} shows that this is true for the classical
wavelet representation on $\ltwor$ (i.e. $n=1$,
$\sigma=\mbox{identity}$, $\theta_1=0$). We will show that the
general case follows from this particular one.
\par
For $M\geq0$, $j<0$, $\psi\in\Psi$, $k\in\mathbb{Z}$ we evaluate
\begin{align*}
|\ip{U^{-M}f}{\tilde{\psi}_{j,k}}|^2&=|\sum_{i=1}^n\int_{\mathbb{R}}U_0^{-M}f_{\sigma^{-M}(i)}(x)N^{j/2}\cj{z_i^kT_0^kU_0^{-j}\psi_{\sigma^{-j}(i)}(x)}\,dx|^2\\
&\leq|\sum_{i=1}^n\int_{\mathbb{R}}|U_0^{-M}f_{\sigma^{-M}(i)}(x)|N^{j/2}|T_0^kU_0^{-j}\psi_{\sigma^{-j}(i)}(x)|\,dx|^2\\
&\leq|\sum_{i=1}^n\int_{\mathbb{R}}U_0^{-M}F_f(x)N^{j/2}T_0^kU_0^{-j}P_{\psi}(x)\,dx|^2,
\end{align*}
where
$$F_f:=\max\{|f_i|\,|\,i\in\{1,...,n\}\},P_{\psi}:=\max\{|\psi_i|\,|\,i\in\{1,...,n\}\},$$
and the last term in the previous inequality is actually
$$=n^2|\ip{U_0^{-M}F_f}{(P_{\psi})^{\tilde{}}_{j,k}}|^2.$$
Apply now lemma 4 of \cite{CSS} for $f\leftrightarrow F_f$ and
$\Psi\leftrightarrow\{P_\psi\,|\,\psi\in\Psi\}$ and with the
previous inequality we obtain the first limit.
\par
The second one can be obtain by a similar argument.
\end{proof}
\par
Theorem \ref{th2_9_0} is one of the crucial steps in obtaining a
characterization for wavelets. The only thing that remains to do
is to use theorem \ref{th2_8} for the quasi-affine system and
carefully compute the Gramians to obtain the next result which
gives the characterizing equation for NTF wavelets in $\ltwor^n$.

\begin{theorem}\label{th2_9}
Let $\Psi=\{\psi_1,...,\psi_p\}$ be a finite subset of $\ltwor^n$.
The following affirmations are equivalent:
\begin{enumerate}
\item
The affine system $X(\Psi)$ is a NTF for $\ltwor^n$;
\item
The following equation hold for almost every $\xi\in\mathbb{R}$,
$i,j\in\{1,...,n\}$:
\\
For all indices $i,j$ with $z_i=z_j$:
\begin{equation}\label{eq2_9_1}
\sum_{\psi\in\Psi}\sum_{m=-\infty}^{\infty}\fo{\psi}_{\sigma^m(i)}(N^m\xi)\cj{\fo{\psi}}_{\sigma^m(j)}(N^m\xi)=\delta_{ij}.
\end{equation}
For all $i,j$ and for all $s\in\mathbb{Z}\setminus N\mathbb{Z}$:
\begin{equation}\label{eq2_9_1_1}
\sum_{\psi\in\Psi}\sum_{m\geq0}\fo{\psi}_{\sigma^{m}(i)}(N^m\xi)\cj{\fo{\psi}}_{\sigma^m(j)}(N^m(\xi+N(\theta_{\sigma^{-1}(i)}-\theta_{\sigma^{-1}(j)})+2s\pi))=0.
\end{equation}
\end{enumerate}
\end{theorem}
\begin{proof}
According to theorem \ref{th2_9_0}, the affine system $X(\Psi)$ is
a normalized tight frame if and only if the quasi-affine
$X^q(\Psi)$ is. This can be reformulated as
$$\{\tilde{\psi}_{m,0}\,|\,m\leq0,\psi\in\Psi\}\cup\{\tilde{\psi}_{m,r}\,|\,m>0,r\in\{0,...,N^{m}-1\},\psi\in\Psi\}$$
is a NTF generator for $\ltwor^n$. We will use theorem \ref{th2_8}
for this family.
\par
We have to compute the Fourier transforms:
$$\fo{\tilde{\psi}}_{m,0}(\xi)=\left(\fo{\psi}_{\sigma^{-m}(i)}(N^{-m}\xi)\right)_{i\in\{1,...,n\}},\quad(m\leq0),$$
and
$$\fo{\tilde{\psi}}_{m,r}(\xi)=N^{-m/2}\left(e^{-i(N^{-m}\xi+\theta_{\sigma^{-m}(i)})r}\fo{\psi}_{\sigma^{-m}(i)}(N^{-m}\xi)\right)_{i\in\{1,...,n\}},$$
when $m>0,r\in\{0,...,N^{m}-1\}.$ Change $m$ into $-m$ and the
equations of theorem \ref{th2_8} are
\begin{eqnarray}
\sum_{\psi\in\Psi}\sum_{m\geq0}\fo{\psi}_{\sigma^m(i)}(N^m(\xi-\theta_i))\cj{\fo{\psi}}_{\sigma^m(i)}(N^m(\xi-\theta_j+2k\pi))+\nonumber\\
+\sum_{\psi\in\Psi}\sum_{m<0}N^m\left(\sum_{r=0}^{N^{-m}-1}e^{-ir(\theta_{\sigma^m(i)}+N^m(\xi-\theta_i)-(\theta_{\sigma^m(j)}+N^m(\xi-\theta_j+2k\pi)))}\cdot\right.\nonumber\\
\left.\cdot
\fo{\psi}_{\sigma^m(i)}(N^m(\xi-\theta_i))\cj{\fo{\psi}}_{\sigma^{m}(j)}(N^m(\xi-\theta_j+2k\pi))\right)=\delta_{ij}\delta_l,\label{eq2_9_3}
\end{eqnarray}
$i,j\in\{1,...,n\},k\in\mathbb{Z}$ for almost all
$\xi\in\mathbb{R}$.
\par
The rest of the proof will be just the attempt to rewrite the
equation (\ref{eq2_9_3}) in the nicer form given in the statement
of the theorem.
\par
We evaluate the inner sum, for fixed $\xi,i,j,k$. Denote by
$$\alpha(m):=\theta_{\sigma^m(i)}-\theta_{\sigma^m(j)}-N^m(\theta_i-\theta_j+2k\pi),\quad(m\leq0).$$
\begin{align*}
S_m&:=N^m\sum_{r=0}^{N^{-m}-1}e^{-ir(\theta_{\sigma^m(i)}+N^m(\xi-\theta_i)-(\theta_{\sigma^m(j)}+N^m(\xi-\theta_j+2k\pi)))}\\
&=\left\{\begin{array}{ccc}
N^m\frac{1-e^{-i\left(N^{-m}(\theta_{\sigma^m(i)}-\theta_{\sigma^m(j)})-(\theta_i-\theta_j)\right)}}{1-\alpha(m)}&\mbox{if}&e^{-i\alpha(m)}\neq
1\\
1&\mbox{if}&e^{-i\alpha(m)}=1.
\end{array}
\right.
\end{align*}
But $m\leq0$ so $N^{-m}\theta_{\sigma^{m}(i)}\equiv\theta_i\mod
2\pi$ and similarly for $j$ so that
$$S_m=\left\{\begin{array}{ccc}
0,&\mbox{if}&\alpha(m)\not\equiv0\mod 2\pi\\
1,&\mbox{if}&\alpha(m)\equiv0\mod 2\pi.
\end{array}
\right.$$ Also note that, if for some $m<0$, $\alpha(m)\in
2\pi\mathbb{Z}$ then $\alpha(m+1)\in 2\pi\mathbb{Z}$.
\par
Indeed,
$$2N\pi\mathbb{Z}\ni
N\alpha(m)=N\theta_{\sigma^{m}(i)}-N\theta_{\sigma^{m}(j)}-N^{m+1}(\theta_i-\theta_j+2k\pi)$$
But $$N\theta_{\sigma^m(i)}\equiv\theta_{\sigma^{m+1}(i)}\mod
2\pi$$ and similarly for $j$ therefore we see that $\alpha(m+1)$
is in $2\pi\mathbb{Z}$ also.
\par
Thus, if for some $m_0\leq 0$ we have
$\alpha(m_0)\in2\pi\mathbb{Z}$ then all $\alpha(m)$ are in
$2\pi\mathbb{Z}$ for $0\geq m\geq m_0$, and so $S_m=1$ for $m\geq
m_0$.
\par
Define
$$m_0:=\min\{m\leq 0\,|\,\alpha(m)\in
2\pi\mathbb{Z}\}\in\{-\infty,...,-1,0\}.$$ Rewrite
(\ref{eq2_9_3}):
\begin{equation}\label{eq2_9_4}
\sum_{\psi\in\Psi}\sum_{m\geq
m_0}\fo{\psi}_{\sigma^m(i)}(N^m(\xi-\theta_i))\cj{\fo{\psi}}_{\sigma^m(j)}(N^m(\xi-\theta_j+2k\pi))=\delta_{ij}\delta_k.
\end{equation}
Of course, $m_0$ depends on $i,j,k$ and we have to make this
dependence more precise.
\par
If $z_i=z_j$ and $k=0$ then
$\theta_{\sigma^{m}(i)}=\theta_{\sigma^{m}(j)}$ for all
$m\in\mathbb{Z}$ so $\alpha(m)=0$ for all $m\leq0$ and therefore
$m_0=-\infty$. In this case the equation (\ref{eq2_9_4}) is
$$\sum_{m=-\infty}^{\infty}\fo{\psi}_{\sigma^m(i)}(N^m(\xi-\theta_i))\cj{\fo{\psi}}_{\sigma^m(j)}(N^m(\xi-\theta_j))=\delta_{ij},$$
which after the change of variable $\xi-\theta_i\mapsto\xi$
becomes (\ref{eq2_9_1}).
\par
If $z_i\neq z_j$ or $k\neq0$ then $m_0$ is finite. When $z_i\neq
z_j$ we have $\theta_{\sigma^{m}(i)}\neq\theta_{\sigma^{m}(j)}$
for all $m$, by cyclicity
$\theta_{\sigma^{m}(i)}-\theta_{\sigma^{m}(j)}$ will assume only a
finite number of nonzero values  in $(-2\pi,2\pi)$ as $m$ varies.
For $m$ close to $-\infty$ the quantity
$-N^m(\theta_i-\theta_j+2k\pi)$ is too small to cover the distance
from $\theta_{\sigma^{m}(i)}-\theta_{\sigma^{m}(j)}$ to any point
in $2\pi\mathbb{Z}$ so that $\alpha(m)\not\in 2\pi\mathbb{Z}$.
\par
We change the variable in (\ref{eq2_9_4}),
$\xi'-\theta_{\sigma^{m_0}(i)}=N^{m_0}(\xi-\theta_i)$. Then
$$N^{m_0}(\xi-\theta_j+2k\pi)=\xi'-\theta_{\sigma^{m_0}(j)}-\alpha(m_0),$$
so (\ref{eq2_9_4}) is equivalent to
\begin{equation}\label{eq2_9_5}
\sum_{m\geq0}\fo{\psi}_{\sigma^{m+m_0}(i)}(N^m(\xi-\theta_{\sigma^{m_0}(i)}))\cj{\fo{\psi}}(N^m(\xi-\theta_{\sigma^{m_0}(j)}-\alpha(m_0)))=0,
\end{equation}
\par
We prove that (\ref{eq2_9_5}) is equivalent to:
\par
For all $i,j$ and $s\in\mathbb{Z}\setminus N\mathbb{Z}$:
\begin{equation}\label{eq2_9_2}
\sum_{\psi\in\Psi}\sum_{m\geq0}\fo{\psi}_{\sigma^{m}(i)}(N^m(\xi-\theta_i))\cj{\fo{\psi}}_{\sigma^m(j)}(N^m(\xi-\theta_j+A(i,j)+2s\pi))=0,
\end{equation}
where
$A(i,j)=N(\theta_{\sigma^{-1}(i)}-\theta_{\sigma^{-1}(j)})-(\theta_i-\theta_j)$.
\par
First, assume (\ref{eq2_9_2}) holds. Fix $i,j\in\{1,...,n\}$ ,
$k\in\mathbb{Z}$ with $z_i\neq z_j$ or $k\neq0$. We prove that
$-\alpha(m_0)$ is of the form
$A(\sigma^{m_0}(i),\sigma^{m_0}(j))-2s\pi$ with
$s\in\mathbb{Z}\setminus N\mathbb{Z}$.
\par
We have
$$\alpha(m_0)=\theta_{\sigma^{m_0}(i)}-\theta_{\sigma^{m_0}(j)}-N^{m_0}(\theta_i-\theta_j+2k\pi).$$
By the definition of $m_0$, $\alpha(m_0-1)\not\in2\pi\mathbb{Z}$.
But observe that
\begin{align*}
N\alpha(m_0-1)&=N(\theta_{\sigma^{m_0-1}(i)}-\theta_{\sigma^{m_0-1}(j)})-N^{m_0}(\theta_i-\theta_j+2k\pi)\\
&=\alpha(m_0)+A(\sigma^{m_0}(i),\sigma^{m_0}(j)).
\end{align*}
Also
$A(\sigma^{m_0}(i),\sigma^{m_0}(j))\in 2\pi\mathbb{Z}$ and
$\alpha(m_0)\in 2\pi\mathbb{Z}$ so $N\alpha(m_0-1)\in
2\pi\mathbb{Z}$ which means that
$\alpha(m_0-1)\in\frac{2\pi}{N}\mathbb{Z}$.
\par
We also know that $\alpha(m_0-1)\not\in 2\pi\mathbb{Z}$, hence
$\alpha(m_0-1)=2\pi\frac{s}{N}$ with $s\in\mathbb{Z}\setminus
N\mathbb{Z}$. Then
$$-\alpha(m_0)=-2s\pi+A(\sigma^{m_0}(i),\sigma^{m_0}(j)).$$
Now using (\ref{eq2_9_2}) for $\sigma^{m_0}(i),\sigma^{m_0}(j)$
and $-s$, we obtain (\ref{eq2_9_5}).
\par
For the converse, suppose (\ref{eq2_9_5}) holds and fix $i,j$ and
$s\in\mathbb{Z}\setminus N\mathbb{Z}$. Let $k:=2s\pi+A(i,j)$ and
we show that the $m_0$ associated to $i,j$ and $k$ is $0$.
\par
Indeed, we compute: $\alpha(-1)=-2\frac{s}{N}\pi\not\in
2\pi\mathbb{Z}$. Since $\alpha(0)=-2s\pi-A(i,j)$, using
(\ref{eq2_9_5}), the relation (\ref{eq2_9_2}) is obtained.
\par
Equation (\ref{eq2_9_2}) is clearly equivalent to
(\ref{eq2_9_1_1}) after a change of variable.
\end{proof}

\begin{remark}
Note that there is some redundancy in the equations of theorem
\ref{th2_9}. Indeed, the relations (\ref{eq2_9_1}) and
(\ref{eq2_9_1_1}) for $i>j$ follow from the corresponding
relations for $i\leq j$ after a conjugation and a change of
variable.
\par
Also, if $i_1,...,i_l$ form a cycle for $\sigma$ then, once we
have equation (\ref{eq2_9_1}) for $i=j=i_1$, the equation
(\ref{eq2_9_1}) will be also true for $i_2,...,i_l$, because one
needs only to change the variable $\xi\leftrightarrow N^p\xi$ for
some appropriate $p$.
\end{remark}
We can apply theorem \ref{th2_9} to the particular case of the
amplification of the classic representation on $\ltwor$. Our
result complements the equations obtained for orthogonal
super-wavelets in \cite{HL}, theorem 5.13.
\begin{corollary}\label{cor2_9_1}
Let $\sigma$ be the identity permutation and $z_1=...=z_n=1$. Then
$\Psi$ is a NTF wavelet if and only if the following equations are
satisfied for almost all $\xi\in\mathbb{R}$:
$$\sum_{\psi\in\Psi}\sum_{m=-\infty}^{\infty}\fo\psi_i(N^m\xi)\cj{\fo{\psi}}_j(N^m\xi)=\delta_{i,j},\quad(i,j\in\{1,...,n\});$$
$$\sum_{\psi\in\Psi}\sum_{m=0}^{\infty}\fo\psi_i(N^m\xi)\cj{\fo{\psi}}_j(N^m(\xi+2s\pi))=0,\quad(i,j\in\{1,...,n\},s\in\mathbb{Z}\setminus N\mathbb{Z}).$$
\end{corollary}
\par
In \cite{HL}, two normalized tight frames $(e_i)_{i\in I}$ for $H$
and $(f_i)_{i\in I}$ for $K$ are called strongly disjoint if their
direct sum $(e_i\oplus f_i)_{i\in I}$ is a normalized tight frame
for $H\oplus K$. Two normalized tight frame wavelets
$\Psi=\{\psi_1,...,\psi_L\}$, $\Psi'=\{\psi'_1,...,\psi'_L\}$ for
$\ltwor$ (note that the cardinality is the same) are strongly
disjoint if their affine systems $X(\psi)$ and $X(\Psi')$ are
strongly disjoint.
\par
If we look at corollary \ref{cor2_9_1}, we see that the fact that
each component $\{\psi_i\,|\,\psi\in\Psi\}$, $(i\in\{1,...,n\})$
is a NTF wavelet for $\ltwor$ is equivalent to the equations with
$i=j$ (and we reobtain in fact the well known characterization of
wavelets in $\ltwor$, see \cite{Bo2}, \cite{Ca} or \cite{HW}).
Therefore the disjointness part is covered by the equations with
$i\neq j$, hence we can use the corollary \ref{cor2_9_1} to
produce a characterization of disjointness:
\begin{corollary}\label{cor2_9_2}
Two NTF wavelets $\Psi:=\{\psi^1,...,\psi^L\}$,
$\Psi':=\{\psi'^1,...,\psi'^L\}$ for $\ltwor$ are strongly
disjoint if and only if
$$\sum_{l=1}^L\sum_{m=-\infty}^{\infty}\fo{\psi^l}(N^m\xi)\cj{\fo{\psi'^l}}(N^m\xi)=0,$$
$$\sum_{l=1}^L\sum_{m=0}^{\infty}\fo{\psi^l}(N^m\xi)\cj{\fo{\psi'^l}}(N^m(\xi+2s\pi))=0,\quad(s\in\mathbb{Z}\setminus N\mathbb{Z}).$$
\end{corollary}
\par
In the remainder of this section we present another application of
the characterization theorem. The following oversampling result is
proved in \cite{CS} for the scaling factor $N=2$: if one has a NTF
wavelet $\psi$ (for $\ltwor$), then also
$\eta:=\frac{1}{p}\psi\left(\frac{x}{p}\right)$ is a NTF wavelet,
where $p$ is any odd number.
\par
We will refine this result here: not only that $\eta$ is a NTF
wavelet, but it is part of a super-wavelet $(\eta,...,\eta)$ for
$\ltwor^p$, and the starting wavelet is an orthogonal one if and
only if the super-wavelet is also orthogonal. Here is the precise
formulation of this result:
\par
Let $p$ be a positive integer which is prime with $N$. For
$\psi\in\ltwor$, define
$$\eta(\psi)(x)=\frac{1}{p}\psi\left(\frac{x}{p}\right),\quad(x\in\mathbb{R}),$$
$$\vec{\eta}(\psi)=(\eta(\psi),...,\eta(\psi))\in\ltwor^p,\quad(\psi\in\Psi).$$
Let $\rho:=e^{-2\pi i/p}$, $z_k=\rho^k$, $k\in\{1,...,p\}$ and
$\sigma$ the permutation of $\{1,...,p\}$ defined by
$(\rho^i)^N=\rho^{\sigma(i)}$ for all $i\in\{1,...,p\}$, and let
$\mathfrak{R}_{\sigma,Z}$ be the corresponding affine structure on
$\ltwor^p$.
\begin{theorem}\label{th2_9_3}
Let $\Psi:=\{\psi_1,...,\psi_p\}$ in $\ltwor$. $\Psi$ is a NTF
wavelet for $\ltwor$ if and only if
$$\vec{\eta}(\Psi):=\{\vec{\eta}(\psi)\,|\,\psi\in\Psi\}$$ is a
NTF wavelet on $\ltwor^p$ for the affine structure
$\mathfrak{R}_{\sigma,Z}$.
\par
Moreover, $\Psi$ is an orthogonal wavelet for $\ltwor$ if and only
if  $\vec{\eta}(\Psi)$ is an orthogonal wavelet for $\ltwor^p$.
\end{theorem}
\begin{proof}
We have to check that the equations of theorem \ref{th2_9} are
satisfied.
\par
The Fourier transform of $\eta(\psi)$ is
$$\fo{\eta(\psi)}(\xi)=\fo{\psi}(px),\quad(\xi\in\mathbb{R}).$$
Then, we can rewrite the equations (\ref{eq2_9_1}) and
(\ref{eq2_9_1_1}) for $\vec{\eta}(\psi)$:
$$\sum_{\psi\in\Psi}\sum_{m=-\infty}^{\infty}|\fo{\psi}|^2(N^mp\xi)=1,$$
$$\sum_{\psi\in\Psi}\sum_{m\geq 0}\fo{\psi}(N^mp\xi)\cj{\fo{\psi}}(N^m(p\xi+a))=0,$$
for all
$a\in\{p(N\theta_{\sigma^{-1}(i)}-N\theta_{\sigma^{-1}(j)}+2s\pi)\,|\,i,j\in\{1,...,p\},s\in\mathbb{Z}\setminus
N\mathbb{Z}\}=:A$. We want to see what this set is and we will
show that it is equal to $\{2q\pi\,|\,q\in\mathbb{Z}\setminus
N\mathbb{Z}\}$.
\par
Note that $\theta_i$ takes all the values of the form $2k\pi/p$
inside $[-\pi,\pi)$ with $k$ integer. Then $p\theta_i$ takes all
values of the form $2k\pi$ inside $[-p\pi,p\pi)$ and therefore
$p\theta_{\sigma^{-1}(i)}-p\theta_{\sigma^{-1}(j)}$ will cover all
values of the form $2k\pi$ within $(-2p\pi,2p\pi)$. Therefore
$$A=\{2\pi(Nl+ps)\,|\,l\in\{-p+1-p+2,...,-1,0,1,...,p-1\},s\in\mathbb{Z}\setminus
N\mathbb{Z}\}.$$
\par
Since $p$ is prime with $N$ it is clear that all numbers in $A$
are of the form $2\pi q$ with $q\in\mathbb{Z}\setminus
N\mathbb{Z}$.
\par
For the converse, take $q$ not divisible by $N$. Since $N$ and $p$
are prime we can write $q=Nl_1+ps_1$ with $l_1,s_1$ integers.
Clearly $s_1$ cannot be divisible by $N$. Also we can write
$l_1=pr_1+l$ with $r_1\in\mathbb{Z}$ and $l\in\{0,...,p-1\}$. So
$$q=Nl+p(Nr_1+s_1)$$
with $l\in\{0,...,p-1\}$, $Nr_1+s_1$ not divisible by $N$ and
therefore $2\pi q$ is in $A$.
\par
Finally, if we change the variable $p\xi\leftrightarrow \xi$, we
see that the equations are equivalent to the fact that $\Psi$ is a
NTF wavelet for $\ltwor$.
\par
The orthogonal case follows by an inspection of the norms: a
normalized tight frame is an orthonormal basis if and only if all
the norms of the vectors are 1. But observe that
$$\|\vec{\eta}(\psi)\|^2=p\|\eta(\psi)\|^2=\|\psi\|^2,$$
and therefore the last equivalence is clear.
\end{proof}

\section{\label{spec}The spectral function and the dimension function}
\begin{definition}\label{def2_10}
A finite subset $\Psi=\{\psi^1,...,\psi^L\}$ of $\ltwor^n$ is
called a semi-orthogonal wavelet if the affine system
$$\{U^jT^k\psi\,|\, j\in\mathbb{Z}, k\in\mathbb{Z} ,\psi\in\Psi\}$$
is a NTF for $\ltworn$ and $W_i\perp W_j$ for $i\neq j$, where
$$W_j=V_{j+1}\ominus V_j=\overline{\operatorname*{span}}\{U^{-j}T^k\psi\,|\,k\in\mathbb{Z},\psi\in\Psi\}=U^{-j}(S(\Psi)),\quad(j\in\mathbb{Z}).$$
and $(V_j)_j$ is the GMRA associated to $\Psi$.
\par
If $\Phi$ is a NTF generator for $V_0$ then $\Phi$ is called a set
of scaling functions for $\Psi$.
\end{definition}

\par
For a semi-orthogonal wavelet, theorem \ref{th2_9_0} and the
orthogonality conditions given in the definition imply the
following:
\begin{proposition}\label{prop2_11}
Let $\Psi$ be a semi-orthogonal wavelet and $(V_j)_j$ its GMRA.
Then the vectors
$$\{\tilde{\psi}_{m,0}=N^{-m/2}U^{m}\psi\,|\,m>0,\psi\in\Psi\}$$
form a NTF generator for $V_0$.
\end{proposition}
\par
This proposition will be one of the essential ingredients for our
results. In $V_0$ we have two NTF generators: one given by the
wavelet $\Psi$ as in proposition \ref{prop2_11} and another given
by a set of scaling functions for $V_0$ which always exists
according to theorem \ref{th2_6}. Computing the local trace for
various operators in two ways will give rise to some fundamental
equalities about wavelets.
\par
The spectral function was introduced for the classical case on
$\ltwor$ by M. Bownik and Z. Rzeszotnik in \cite{BoRz}. It
contains a lot of information about the shift invariant subspace
and it has several nice features. In some sense, the spectral
function corresponds to the diagonal entries of the dual Gramian.
\begin{definition}\label{def2_12}
Let $V$ be a shift invariant subspace of $\ltwor^n$. For
$i\in\{1,...,n\}$, the $i$-th spectral function of $V$ is
$$\sigma_{V,i}(\xi):=\tau_{V,P_{0,i}}(\xi),\quad(\xi\in\mathbb{R}),$$
where $P_{0,i}$ is the projection onto the $0$-th component of the
$i$-th vector
$$(P_{0,i}v)(k,j)=v_{0,i}\delta_{0k}\delta_{ij},\quad\mbox{ for
all }j\in\{1,...,n\}, k\in\mathbb{Z}.$$
\end{definition}
\par
The next proposition is just a consequence of theorem \ref{th3_3}
\begin{proposition}\label{prop2_13}
Let $V$ be a shift invariant subspace and $\Phi$ be a NTF
generator for $V$. Then
$$\sigma_{V,i}(\xi)=\sum_{\varphi\in\Phi}|\fo{\varphi_i}|^2(\xi-\theta_i),\quad(\xi\in\mathbb{R},i\in\{1,...,n\}).$$
\end{proposition}
\par
Using the NTF generator given in proposition \ref{prop2_11}, we
obtain a formula for the spectral functions associated to
wavelets:
\begin{proposition}\label{prop2_14}
Let $\Psi$ be a semi-orthogonal wavelet and $(V_j)_j$ its GMRA.
Then the spectral functions of the core space $V_0$ can be
computed by:
$$\sigma_{V_0,i}(\xi)=\sum_{\psi\in\Psi}\sum_{m\geq1}|\fo{\psi}_{\sigma^m(i)}|^2(N^m(\xi-\theta_i)),\quad(\xi\in\mathbb{R},i\in\{1,...,n\}).$$
\end{proposition}
Combining proposition \ref{prop2_13} and \ref{prop2_14} we obtain
a Gripenberg-Weiss-type formula:
\begin{corollary}\label{cor2_15}
If $\Psi$ is a semi-orthogonal wavelet and $\Phi$ is a set of
scaling function for it, then
$$\sum_{\varphi\in\Phi}|\fo{\varphi_i}|^2(\xi-\theta_i)=\sum_{\psi\in\Psi}\sum_{m\geq1}|\fo{\psi}_{\sigma^m(i)}|^2(N^m(\xi-\theta_i)),\quad(\xi\in\mathbb{R},i\in\{1,...,n\}).$$
\end{corollary}
\par
The dimension function was extensively used on $\ltwor$ to decide
whether a wavelet is a MRA wavelet or not, see \cite{HW}. It also
gives information about the multiplicity of the wavelets
(\cite{B},\cite{BM},\cite{BMM},\cite{Web}).

\begin{definition}\label{def2_16}
Let $V$ be a shift invariant subspace of $\ltwor^n$ and $J_{per}$
its range function. The dimension function of $V$ is
$$\mbox{dim}_V(\xi)=\mbox{dim}(J_{per}(\xi)),\quad(\xi\in\mathbb{R}).$$
\end{definition}
Since the trace of a projection is the dimension of its range,
using the definition of the local trace we have the following
equality:
\begin{proposition}\label{prop2_17}
Let $V$ be a shift invariant subspace of $\ltwor^n$. Then
$$\mbox{dim}_V(\xi)=\tau_{V,I}(\xi),$$
where $I$ is the identity operator on $\ltwoz^n$.
\end{proposition}
Then, with theorem \ref{th3_3}, and proposition \ref{prop2_13} we
have
\begin{proposition}\label{prop2_18}
If $\Phi$ is a NTF generator for a shift invariant subspace $V$
then
$$\mbox{dim}_V(\xi)=\sum_{\varphi\in\Phi}\sum_{i=1}^n\sum_{k\in\mathbb{Z}}|\fo{\varphi}_i|^2(\xi-\theta_i+2k\pi)=
\sum_{i=1}^n\operatorname*{Per}(\sigma_{V,i})(\xi),\quad(\xi\in\mathbb{R}).$$
where
$$\operatorname*{Per}(f)(\xi)=\sum_{k\in\mathbb{Z}}f(\xi+2k\pi),\quad(\xi\in\mathbb{R},f\in\loner).$$
\end{proposition}
We can choose the scaling functions to be as in theorem
\ref{th2_6}, and the consequence is that the dimension function
and the multiplicity function are the same (a generalization of
the result of E. Weber \cite{Web}).
\begin{theorem}\label{th2_19}
If $\Psi$ is a semi-orthogonal wavelet and $V_0$ the core space of
the associated GMRA. Then the multiplicity function and the
dimension function of $V_0$ are equal and they are given by
$$m_{V_0}(\xi)=\mbox{dim}_{V_0}(\xi)=\sum_{\psi\in\Psi}\sum_{m\geq1}\sum_{i=1}^n\sum_{k\in\mathbb{Z}}|\fo{\psi}_{\sigma^m(i)}|^2(N^m(\xi-\theta_i+2k\pi)=:D_{\Psi}(\xi),$$
for almost all $\xi\in [-\pi,\pi)$.
\end{theorem}
\par
Next we want to establish some dilation formulas for the spectral
function and for the dimension function. We use theorem
\ref{th2_7}.
\begin{proposition}\label{prop2_20}
Let $V$ be a SI subspace of $\ltwor^n$. Then for all
$j\in\{1,...,n\}$,
$$\sigma_{U^{-1}V,j}(\xi)=\sigma_{V,\sigma^{-1}(j)}\left(\frac{\xi-\theta_j+N\theta_{\sigma^{-1}(j)}}{N}\right),\quad(\xi\in\mathbb{R}),$$

$$\mbox{dim}_{U^{-1}V}(\xi)=\sum_{l=0}^{N-1}\mbox{dim}_V\left(\frac{\xi+2l\pi}{N}\right),\quad(\xi\in\mathbb{R}).$$
\end{proposition}

\begin{proof}
With theorem \ref{th2_7}, we have to compute for
$l\in\{0,...,N-1\}$, $\mathcal{S}^*D_l^*P_{0,j}D_l\mathcal{S}$ is
the rank-one operator given by the vector
$\mathcal{S}^*D_l^*e_{0,j}$. We have for
$k\in\mathbb{Z},i\in\{0,...,N-1\}$:
$$\mathcal{S}^*D_l^*e_{0,j}(k,i)=(D_l^*e_{0,j})(k,\sigma(i))=e_{0,j}(2l\pi+\theta_{\sigma(i)}-N\theta_i+2kN\pi),$$
and this shows that
$$\mathcal{S}^*D_l^*e_{0,j}=\left\{\begin{array}{ccc}
e_{k_j,\sigma^{-1}(j)}&\mbox{if}&l=l_j\\
0,& &\mbox{otherwise},
\end{array}\right.$$
where $l_j\in\{0,...,N-1\}$ and $k_j\in\mathbb{Z}$ are the only
such numbers that verify the equation
$2l_j\pi+\theta_j-N\theta_{\sigma^{-1}(j)}+2k_jN\pi=0$. Therefore
$$\mathcal{S}^*D_l^*P_{0,j}D_l\mathcal{S}=\left\{\begin{array}{ccc}
P_{k_j,\sigma^{-1}(j)}=\lambda(k_j)P_{0,\sigma^{-1}(j)}\lambda(k_j)^*&\mbox{if}&l=l_j,\\
0,& &\mbox{otherwise}. \end{array}\right.$$ Using theorem
\ref{th2_7} and then the periodicity property given in proposition
\ref{prop3_4}, we obtain:
\begin{align*}
\sigma_{U^{-1}V,j}(\xi)&=\tau_{V,P_{k_j,\sigma^{-1}(j)}}\left(\frac{\xi+2l_j\pi}{N}\right)\\
&=\sigma_{V_0,\sigma^{-1}(j)}\left(\frac{\xi+2l_j\pi}{N}+2k_j\pi\right),
\end{align*}
and, with the definition of $k_j$ and $l_j$ the equation follows.
\par
The equation for the dimension functions follows from the fact
that $\mathcal{S}^*D_l^*ID_l\mathcal{S}=I$ for all
$l\in\{0,...,N-1\}$.
\end{proof}

\begin{proposition}\label{prop2_21}
Let $(V_j)_j$ be a GMRA and $\Phi$ a NTF generator for $V_0$.
Then, for almost every $\xi\in\mathbb{R}$ and all
$i\in\{1,...,n\}$,
$$\lim_{m\rightarrow\infty}\sum_{\varphi\in\Phi}|\fo{\varphi}_i|^2\left(\frac{\xi}{N^m}\right)=1.$$
\end{proposition}
\begin{proof}
With the monotone convergence theorem \ref{th4_3}, we have that
$\sigma_{U^{-m}V,i}(\xi)$ converges pointwise a.e. to
$\sigma_{\ltwor^n,i}(\xi)=1$. With proposition \ref{prop2_20}, by
induction we get:
\begin{align*}
\sigma_{U^{-m}V,i}(\xi)&=\sigma_{V_0,\sigma^{-m}(i)}\left(\frac{\xi-\theta_i+N^m\theta_{\sigma^{-m}(i)}}{N^m}\right)\\
&=\sum_{\varphi\in\Phi}|\fo{\varphi}_{\sigma^{-m}(i)}|^2\left(\frac{\xi-\theta_i}{N^m}\right).
\end{align*}
Changing $\xi-\theta_i$ to $\xi$ we obtain
$$\lim_{m\rightarrow\infty}\sum_{\varphi\in\Phi}|\fo{\varphi}_{\sigma^{-m}(i)}|^2\left(\frac{\xi}{N^m}\right)=1.$$
But $\sigma$ is a finite permutation so $\sigma^{p}=\mbox{id}$ for
some $p$. Thus, for the subsequence $mp$, $\sigma^{-mp}(i)=i$ for
all $m$, hence
$$\lim_{m\rightarrow\infty}\sum_{\varphi\in\Phi}|\fo{\varphi}_{i}|^2\left(\frac{\xi}{N^{mp}}\right)=1,$$
and if we apply this to $\xi$, $\xi/N$,...,$\xi/N^{p-1}$ we obtain
the desired limit.
\end{proof}
\par
The limit given in proposition \ref{prop2_21} enables us to give a
lower bound for the dimension function of a wavelet. We can
conclude that, if there are two distinct indices $i,j$ such that
$z_i=z_j$ then the dimension function of any wavelet is strictly
bigger then $1$ at some points so there are no MRA wavelets, and
for any wavelet one needs at least 2 scaling functions. If there
are three indices for which the points $z_i$ are the same, then
any wavelet needs at least three scaling functions, and so on.
This result generalizes and refines proposition 5.16 in \cite{HL}.
\begin{corollary}\label{cor2_22}
If $\Psi$ is a semi-orthogonal wavelet then, for all $\alpha_0\in
[-\pi,\pi)$,
$$\limsup_{\xi\rightarrow\alpha_0}D_{\Psi}(\xi)\geq\mbox{card}\{i\in\{1,...,n\}\,|\,\theta_i=\alpha_0\}.$$
In particular, if there are two distinct indices $i\neq j$ with
$\theta_i=\theta_j$ then there are no MRA wavelets.
\end{corollary}
\begin{proof}
Let $\Phi$ be a set of scaling functions associated to $\Psi$ and
compute $D_{\Psi}$ at $\xi/N^m+\alpha_0$ using proposition
\ref{prop2_18}. Then take the limsup and use proposition
\ref{prop2_21}.
\par
For a MRA wavelet $D_{\Psi}$ is constant $1$, and the second
statement follows by contradiction.
\end{proof}

\end{document}